\documentclass{amsart}

\usepackage{latexsym}
\usepackage{amssymb}
\usepackage{graphicx}

\newtheorem{THM}{Theorem}

\newcommand{\bR}{\mathbb{R}}
\newcommand{\be}{\begin{equation}}
\newcommand{\ee}{\end{equation}}
\newcommand{\mS}{\mathbb{S}}

\newcommand{\Zed}{\mathbb{Z}}
\newcommand{\ol}[1]{\overline{#1}}

\newcommand{\lp}{\left(}
\newcommand{\rp}{\right)}

\newcommand{\lc}{\left\{}
\newcommand{\rc}{\right\}}
\newcommand{\lab}{\left|}
\newcommand{\rab}{\right|}
\newcommand{\Lap}{\Delta}

\newcommand{\Prob}{\mathbb{P}}

\DeclareMathOperator{\sgn}{sign}

\newcommand{\Ind}{\boldsymbol 1}

\newcommand{\Iso}{\Zed/2\Zed\times\mathrm{SO}(2)}

\newcommand{\R}{{\mathbf R}}
\newcommand\Tu{\mathbb{T}}
\newcommand{\ba}[1]{\begin{array}{#1}}
\newcommand{\ea}{\end{array}}
\newcommand{\al}{\alpha}
\newcommand{\metr}{g}
\newcommand{\mcil}{M_{\text{cylinder}}}
\newcommand{\mcono}{M_{\text{cone}}}
\newcommand{\misura}{\omega}
\newcommand{\sing}{\mathcal Z}
\newcommand{\cc}[1]{{C^\infty_c\left(#1\right)}}
\newcommand{\inm}{M}
\newcommand{\delcomp}{\Delta|_{\cc \inm}}
\newcommand{\ldue}{L^2(M,\misura)}

\newcommand{\sobduec}{H^2\left(M_\alpha,\misura\right)}

\newcommand{\deld}{\delf}
\newcommand{\delf}{\Delta_F}

\begin{document}

\title[Brownian motion and singularities]{Extensions of Brownian motion to a family of Grushin-type singularities}
\author{Ugo Boscain \and Robert W.\ Neel}
\address{CNRS, Sorbonne Universit\'e, Inria, Universit\'e de Paris, Laboratoire Jacques-Louis Lions, Paris, France}
\email{ugo.boscain@upmc.fr}
\address{Department of Mathematics, Lehigh University, Bethlehem, PA, USA}
\email{robert.neel@lehigh.edu}
\begin{abstract} 

We consider a one-parameter family of Grushin-type singularities on surfaces, and discuss the possible diffusions that extend Brownian motion to the singularity. This gives a quick proof and clear intuition for the fact that heat can only cross the singularity for an intermediate range of the parameter. When crossing is possible and the singularity consists of one point, we give a complete description of these diffusions, and we describe a ``best'' extension, which respects the isometry group of the surface and also realizes the unique symmetric one-point extension of the Brownian motion, in the sense of Chen-Fukushima. This extension, however, does not correspond to the bridging extension, which was introduced by Boscain-Prandi, when they previously considered self-adjoint extensions of the Laplace-Beltrami operator on the Riemannian part for these surfaces. We clarify that several of the extensions they considered induce diffusions that are carried by the Marin compactification at the singularity, which is much larger than the (one-point) metric completion. In the case when the singularity is more than one-point, a complete classification of diffusions extending Brownian motion would be unwieldy. Nonetheless, we again describe a ``best'' extension which respects the isometry group, and in this case, this diffusion corresponds to the bridging extension. A prominent role is played by Bessel processes (of every real dimension) and the classical theory of one-dimensional diffusions and their boundary conditions.
\end{abstract}

\maketitle

\section{Introduction}

Consider the (open) Riemannian manifold $(M,\metr)$ where $M=\lp\R\setminus\{0\}\rp\times\Tu$ (here $\Tu$ is the one-dimensional torus), and
\begin{equation}
\metr=dx^2+|x|^{-2\al} d\theta^2, \mbox{ that is, in matrix notation }\metr=\left(\ba{cc} 1&0\\0&|x|^{-2\al}\ea\right).\label{u-metric}
\end{equation}
Here $x\in\R\setminus\{0\}$, $\theta\in \Tu$ and $\alpha\in\R$ is a parameter. An orthonormal frame for the metric \eqref{u-metric} is given by the pair of vector fields
$$
X=\left(
\begin{array}{c}1\\0\end{array}\right),~~~Y=\left(\begin{array}{c}0\\|x|^\alpha \end{array}\right).
$$
Define 
\begin{gather*}
	\mcil=\R\times\Tu, \qquad \mcono=\mcil/\sim,
\end{gather*}
where $(x_1,\theta_1)\sim(x_2,\theta_2)$ if and only if $x_1=x_2=0$. 

When $\al\geq 0$, extending the vector fields $X$ and $Y$ to $\mcil$, the natural control-theoretic notion of the length of a curve shows that there are paths of finite length between $M^+=\{x>0\}\times\Tu$ and $M^-=\{x<0\}\times\Tu$ (which are the Riemannian geodesics for $x\neq 0$ and which are tangent to $X$ when $x=0$), and, as explained in \cite{UgoDario},  this extended distance makes $\mcil$  into a metric space (and a length space) in a way that induces on $\mcil$ its original topology. 
Similarly, when $\al< 0$ the distance induced by the metric \eqref{u-metric} extend naturally to $\mcono$. 
Said differently, $\mcil$ (for $\al\geq0$) and  $\mcono$ (for $\al<0$) give the metric compactifications of $M$ with respect to this distance. We denote these metric spaces by $M_{\alpha}$
($M_{\alpha}=\mcil$ if $\alpha\geq0$ and $M_{\alpha}=\mcono$ if $\alpha<0$) and note that $M$ is the Riemannian subset, while we let $\sing=M_{\alpha}\setminus M=\{x=0\}$ be the singular set (terminology which will be further justified in a moment).

This construction gives a one-parameter family of natural singularity models relevant to rank-varying sub-Riemannian/almost-Riemannian geometry that includes the well-known Grushin cylinder (the obvious quotient of the Grushin plane), when $\alpha=1$. Moreover, the case $\alpha\geq1$ corresponds to  an almost-Riemannian structure in the sense of \cite{agrachevboscainsigalotti,boscain-laurent} and of \cite[Chapter 9]{book}. 

However, even though the metric (and length minimizing curves) extends across the singularity,  the Riemannian metric $g$  (except for $\al=0$, which gives a standard cylinder) is singular on $\sing$, and this is also the case for the Riemannian volume $\omega$ and for the Laplace-Beltrami operator, that take the form
\begin{align*}
\misura&=\sqrt{\det g} \,dx\,d\theta=|x|^{-\alpha}dx\,d\theta, \text{ and}\\
\Delta&=\frac{1}{\sqrt{\det g}}\sum_{j,k=1}^2 \partial_j \left(\sqrt{\det g} \, g^{jk}\partial_k\right) =\partial_x^2+|x|^{2\alpha}\partial_\theta^2u-\frac{\alpha}{x}\partial_x .
\end{align*}

Thus, if one wishes to consider the heat equation or Schr\"odinger equation on $M$, or to consider Brownian motion, one must consider the behavior at the boundary. From the perspective of functional analysis, this means considering self-adjoint extensions of the Laplacian on $M$. Indeed, in \cite{UgoDario}, the following two basic results were proven.

\begin{THM}[\cite{UgoDario}]\label{thm:sa}
 The operator $\delcomp$ is essentially self-adjoint in $\ldue$ if and only if $\alpha\in (-\infty,-3]\cup[1,\infty)$.
 \end{THM}
 
An immediate consequence  of Theorem \ref{thm:sa} is that for $\alpha\in (-\infty,-3]\cup[1,\infty)$ the only self-adjoint extension of $\Delta$ is 
the Friedrich extension $\deld$.

If $\alpha\notin(-\infty,-3]\cup[1,\infty)$ the next theorem gives some  additional information.

\begin{THM}[\cite{UgoDario}]\label{thm:sa2}
Let $\widehat \Delta$ be the Fourier transform of $\Delta$ in the variable $\theta$. We have the following
\begin{itemize}
\item if $\alpha\in(-3,-1]$, only the first  Fourier component of $\widehat \Delta$ is not essentially self-adjoint.
\item if $\alpha\in(-1,1)$, all the Fourier components of  $\widehat \Delta$ are not essentially self-adjoint.
\end{itemize}
 \end{THM}
When considering the heat equation $
\partial_t\phi=\Delta\phi
$
in $\ldue$, a consequence  of Theorem \ref{thm:sa2} is that 
\begin{itemize}
\item when $\alpha\in(-3,-1]$ there are self-adjoint extensions of 
$\Delta$ that permit  only the average over $\Tu$ of $\phi$ to flow through $\sing$. However as explained in \cite{UgoDario} the only  Markovian extension of $\Delta$ is $\deld$ 
that does not permit any communication between $M^+$and $M^-$.

 \item when $\alpha\in(-1,1)$ there are self-adjoint extensions of 
$\Delta$ that permit  full communication between $M^+$ and $M^-$ and that are Markovian. In particular, there is a self-adjoint extension called the {\em bridging extension}  realizing the maximal communication between the two sides, the domain of which is
$$
\lc \sobduec \Big\vert u(0^+,\cdot)=u(0^-,\cdot),\, \lim_{x\rightarrow 0^+} |x|^{-\alpha}\partial_x u(x,\cdot)=\lim_{x\rightarrow 0^-} |x|^{-\alpha}\partial_x u(x,\cdot) \rc ,
$$
where $\sobduec=\{u\in L^2(M,\omega),~~|\nabla u|,\Delta u\in     L^2(M,\omega)\}$.
 \end{itemize}

The purpose of the present note is to consider diffusions on $M_{\alpha}$ that extend Brownian motion on $M$. One aspect of this is to give the path properties that correspond to many of the above results. For example, the fact that there are Markov extensions allowing communication between $M^+$ and $M^-$ exactly when $\alpha\in (-1,1)$ corresponds to the fact that for $\alpha<-1$, Brownian motion on $M$ never hits the singularity, and thus cannot cross it, while for $\alpha>1$, is an exit-only boundary for $M$ (essentially in the sense of the classical Feller classification), so the process must be absorbed at the singularity (assuming it is conservative and cannot be killed) and thus also cannot cross it. In particular, the $x$-marginal of Brownian motion on $M$ is given by a Bessel process of dimension $d=1-\alpha$, so both the behavior of the process at the singularity, as well as the stochastic completeness near infinity, follow, and thus stochastic methods provide elementary proofs and intuition for these results. From the other side, this one-parameter family of geometries provides examples ``in nature''  of Bessel processes of all real dimensions (such examples also arise in SLE, but not naturally in Riemannian geometry). This material is treated in Section \ref{Sect:Bessel}.

When considering the possible extensions for $\alpha\in(-1,1)$, it is important to note that the Martin compactification of $M$ at the singularity (in what follows, we are concerned with the behavior at the singularity, and thus all of our compactifications are done there, ignoring what happens near infinity, since there the structure is Riemannian, not singular) is larger than $\sing$, that is, larger than the metric compactification (at the singularity). Various self-adjoint extensions of $\Delta$ mentioned above are carried by the Martin boundary (at the singularity). For example, Neumann boundary conditions make the process undergo instantaneous normal reflection at the singularity, back into the component of $M$ it came from. But such an extension clearly cannot descend to a strong Markov process on $M_{\alpha}$. Thus, here we treat extensions that are carried by $M_{\alpha}$ itself, so that our results differ from, and complement, those of \cite{UgoDario}. We also do not restrict our attention to symmetric extensions. (In this connection, it is worth mentioning that we do not treat non-Markov extensions, so we have no contribution to the above results for $\alpha\in (-3,-1]$.)

In the case $\alpha\in (-1,0)$ when $M_{\alpha}$ is a topological cone, we are able to give a complete description of (conservative) diffusions on $M_{\alpha}$ the extend Brownian motion on $M$. It is worth noting that these correspond to one-point extensions in the sense of Chen and Fukushima (see \cite{ChenFukushimaBook,ChenFukushimaPaper}), and we identify the unique symmetric extension spending zero time at the singularity. For this extension, we see that only the average over $\Tu$ of a function flows through the singularity under the corresponding semigroup (note that in this case, the bridging extension does not correspond to a diffusion on $M_{\alpha}$). Thus the same phenomenon observed in \cite{UgoDario} for non-Markov self-adjoint extensions of $\Delta$ for $\alpha\in(-3,-1]$ is replicated here for Markov processes that respect the topology of $M_{\alpha}$ when $\alpha\in(-1,0)$. This is carried out in Section \ref{Sect:Cone}.

In the case $\alpha\in [0,1)$ when $M$ is a topological cylinder, the larger singular set makes a complete classification of extensions complicated. However, one can describe the basic features, and we also construct the unique diffusion spending time 0 at $\sing$ and respecting the symmetries of $M_{\alpha}$, and show that in this case it corresponds to the bridging extension. This comprises Section \ref{Sect:Cylinder}.

The first author was supported by the ANR project SRGI ANR-15-CE40-0018 and by the ANR project Quaco ANR-17-CE40-0007-01. The second author was partially supported by grant \#524713 from the Simons Foundation and by the National Security Agency under Grant Number H98230-15-1-0171. We thank Dario Prandi, Masha Gordina, and Nate Eldredge for helpful conversations.

\section{Bessel processes, stochastic completeness, and boundary conditions}\label{Sect:Bessel}

From the above description of the metric and the induced Laplacian on $M$, we see that, in the $(x,\theta)$ coordinates, Brownian motion evolves by the system of SDEs
\[\begin{split}
dx_t &= dW_t^1 -\frac{\alpha}{2x}\, dt \\
d\theta_t &= \lab x\rab^{\alpha} \, dW_t^2 ,
\end{split}\]
at least until $T_0$, the first hitting time of $\{x=0\}=\sing$, where $W_t^1$ and $W_t^2$ are independent one-dimensional Brownian motions. (The SDE for $x_t$ should be understood as giving a local semi-martingale on $(0,\infty)$ for general real $\alpha$, but the extension until $T_0$ is standard, say, by squared-Bessel processes as mentioned below, or by explicit construction of the transition density, etc.) It is the $x_t$ process that mainly interests us in this section. Note that its evolution does not depend on $\theta_t$, (except possibly on the singular set) so that the situation reduces to a one-dimensional problem. Moreover, observe that, for $x_t>0$ (equivalently, on $M^{+}$), the SDE satisfied by $x_t$ is just that of a Bessel process of dimension $d= 1-\alpha$ (and for $x_t<0$, it is just $-1$ times such a process).

We could take the perspective that $x=0$ gives an interior singular point of the diffusion, and we will below. However, taking advantage of the reflection symmetry, here we can instead consider the process $z_t=x_t^2$ on $[0,\infty)$, so that the singular set corresponds to an included boundary point. Also, the squared Bessel processes are true semi-martingales for all values of $\alpha$, and satisfy the SDE $dz_t = 2\sqrt{z_t} \, dW_t +\lp 1-\alpha\rp\, dt$ until $T_0$. Further, we note that (the law of) $z_t$ determines (the law of) $x_t$ up to the sign of each excursion of $x_t$ away from 0, that is, away from $\sing$.

The behavior of (squared) Bessel processes is well-understood. First of all, the process doesn't explode to infinity (in finite time) for any value of $\alpha$. For $\alpha\leq -1$, $0$ is an entrance-only boundary (in the standard Feller classification for one-dimensional diffusions). For $\alpha\in (-1,1)$, $0$ is a regular boundary, and thus one needs to specify boundary conditions. For $\alpha\geq 1$, $0$ is an exit-only boundary.

Thus, if we do not allow killing at $\sing$ (or anywhere else), $M_{\alpha}$ is stochastically complete for all $\alpha\in \bR$. (This contrasts slightly with \cite{UgoDario}, since they consider only self-adjoint extensions and thus kill the process at $\sing$ when it is an exit-only boundary, instead of letting it be absorbed. But the real point is that the process never explodes to infinity in finite time.) If $\alpha\leq -1$, then Brownian motion on $M$ (almost surely) never hits $\sing$, so there is no need for an extension. The only caveat is if we wish to start the process from $\sing$. In this case, once the process enters $M$, it never returns to $\sing$, but one still needs to specify an entrance law. This, however, amounts to a simpler version of the $\alpha\in(-1,0)$ case, and we briefly treat it in Section \ref{Sect:Short}.

If $\alpha\in (-1,1)$, then Brownian motion on $M$ almost surely hits $\sing$ in finite time, but is then able to leave and re-enter $M$. Thus extensions of Brownian motion are determined by the boundary behavior at $\sing$ and are not unique. The discussion of this case occupies almost all of the remainder of this note. Finally, if $\alpha\geq 1$, there is only one (conservative) extension of Brownian motion on $M$; it is adsorbed at $\sing$ (which is a single point) at time $T_0$, which is almost surely finite.

This not only recovers the stochastic completeness and difference in (Markov) extensions of $\Delta$ depending on whether $\alpha\in (-1,1)$ or not from \cite{UgoDario}, but also allows other properties of the heat flow on $M_{\alpha}$ to deduced from known properties of Bessel processes. For example, for the Grushin cylinder (or Grushin plane), which corresponds to $\alpha=1$, the rate at which heat is absorbed at $\sing$ is given by the transition measure for a 0-dimensional Bessel process (see Section A.2 of \cite{AnjaYor}, for example).

\section{The case $-1< \alpha < 0$}\label{Sect:Cone}

When $-1 < \alpha < 0$, $M_{\alpha}$ has a cone structure at the singularity, and the singularity reduces to a single point. As we saw above, the singularity is a regular boundary for the $|x_t|$ process, so that any diffusion a.s.\ hits the singular point in finite time, and it is possible for the diffusion to leave the singular point. In this situation, there are many possible diffusions extending Brownian motion on $M$ to all of $M_{\alpha}$, but the singularity is simple enough that we can describe them all.

\subsection{Classification of diffusions on $M_{\alpha}$}
Because the singularity is a single point, the behavior of the diffusion at the singularity doesn't depend on the $\theta_t$ process, which means that $x_t$ is a one-dimensional diffusion. As mentioned, the theory of one-dimensional diffusions is completely understood (see, for example, \cite{Knight}). Thus we can give all possible (conservative) diffusions in this case, and this is the first step in determining the possible diffusions on $M_{\alpha}$. Essentially, the possible $x_t$-diffusions depend on two parameters, the degree of ``stickiness'' at $0$ and the skewness at $0$. We also note that $x_t$ is a Bessel process of dimension between 1 and 2 (at least until it hits 0, at which point we don't necessarily instantaneously reflect it), and thus $x_t$ is a semi-martingale.

More precisely, in the classification scheme of It\^o and McKean \cite{ItoMcKean}, $0$ can be a regular point, a left or right shunt, or a trap. Most interesting for us is when $0$ is a regular point, so that the process can cross $0$ in either direction. In this case, the $x_t$-diffusion is determined by its scale function $s(x)$ and its speed measure $m$. Further, the diffusion must agree with the appropriate Bessel process on $\bR\setminus\{0\}$. Hence the scale function is determined up to affine transformations on each of $\{x<0\}$ and $\{x>0\}$, subject to the additional constraint that it is continuous.  We normalize $s$ by translation so that $s(0)=0$. Starting from the ``standard'' scale function $s(x)=x^{\alpha+1}$ for a Bessel process of dimension $1-\alpha\in(1,2)$ (see Section 11.1 of \cite{RevuzYor}, also for the speed measure of a Bessel process which we are about to use), we see that the most general normalized scale function for $x_t$ with $s(0)=0$ is 
\[
s(x)= \begin{cases}
-a(-x)^{1+\alpha} & \text{for $x<0$} \\
(1-a)x^{1+\alpha} & \text{for $x\geq 0$}
\end{cases}
\]
for $0<a<1$. Here we see that $a$ gives the skewness of $x_t$ at $0$, in the sense that, for any $y>0$,
\[
\Prob\lp\text{$x_t$ starting from 0 hits $y$ before $-y$}\rp = a .
\]

Continuing, the speed measure is uniquely determined on $\bR\setminus\{0\}$ by the speed measure of a Bessel process and the above choice of scaling function (as having density $2/s'$ with respect to Lebesgue measure), so that the most general speed measure for $x_t$ is
\[
m= \frac{2}{a(1+\alpha)(-x)^{\alpha}} \Ind_{\{x<0\}} \,dx +\gamma \delta_0+
\frac{2}{(1-a)(1+\alpha)x^{\alpha}} \Ind_{\{x>0\}} \,dx
\]
for some $\gamma\in[0,\infty)$, where $dx$ denotes Lebesgue measure and $\delta_0$ a point mass at $x=0$. Here $\gamma$ gives the degree of ``stickiness'' at 0, in the sense that if $\gamma=0$, the set $\{t>0: x_t=0\}$ almost surely has Lebesgue measure 0, whereas if $\gamma>0$, this set has positive measure. Note that for $\gamma=0$ and $0<a<1$, $x_t$ will be a skew Brownian motion (one can see the survey \cite{LejaySurvey} for a detailed introduction to skew Brownian motion).

A perhaps more appealing (and slightly more general) way to describe $x_t$ is as follows. While $0$ is a regular point for $0<a<1$, $0$ is a left shunt if $a=0$, and a right shunt if $a=1$. (If $\gamma=\infty$, then $0$ is a trap, viewed as an interior point of $M_{\alpha}$ rather than as a boundary point of $[0,\infty)$ as above.)

Next, we consider the $\theta$-process. Recall that this is an $\Tu$-valued process that satisfies the SDE $d\theta_t=|x|^{\alpha} dW_t^2$ for $x\neq 0$. Suppose $x_0>0$, and make the change of variables $y(x)=\frac{1}{\alpha+1}x^{\alpha+1}$, which puts $x_t$ on its natural scale. Let $T_0$ be the first hitting time of 0 for $x_t$ (and thus the first time the process on $M_{\alpha}$ hits the singularity). Then the process satisfies the system of SDEs
\begin{equation}\label{Eqn:NaturalScale}\begin{split}
dy_t &= (1+\alpha)^{\alpha/(1+\alpha)}y^{\alpha/(1+\alpha)}\, dW_t^1 \\
d\theta_t &= (1+\alpha)^{\alpha/(1+\alpha)} y^{\alpha/(1+\alpha)} \, dW_t^2 ,
\end{split}\end{equation}
on the time interval $[0,T_0]$. Note that $y_t$ is a time-changed Brownian motion that we know a.s.\ hits 0 in finite time, and thus it a.s.\ accumulates finite quadratic variation $\int_0^{T_0} (1+\alpha)^{2\alpha/(1+\alpha)} y^{2\alpha/(1+\alpha)} \, dt$ over $[0,T_0]$. Since the quadratic variation of $\theta_t$ on $[0,T_0]$ is equal to that of $y_t$, it is also a.s.\ finite. Further, $\theta_t$ is a martingale and thus a time-changed Brownian motion, and it follows that $\theta_t$ a.s.\ has a limit as $t\nearrow T_0$. In particular, the existence of a limiting angle (which we can think of as an exit angle, as the process exits $M^+$) implies that the invariant sigma-algebra of this stopped process, and thus also the Martin boundary of $M^+$, is non-trivial, even though the process converges to a single point on $M_{\alpha}$ as $t\nearrow T_0$. (We discuss this further below.) More precisely, this holds when the process is started from any initial point with $x\neq 0$, since $M$ is symmetric under reflection in $x$.

More important, we consider how the process leaves the singularity. Because the law of the Brownian excursion is preserved by time-reversal, it follows from the above that $y_t$ a.s.\ accumulates finite quadratic variation on each excursion from 0, and thus $\theta_t$ does as well. Hence, starting from the singularity, $\theta_t$ must have a limit as $t\searrow 0$, that is, the process leaves the singularity, and enters $M=M\setminus\{x=0\}$, with an entrance angle. Thus, let $\mu^+$ and $\mu^-$ be two probability measures on $\Tu$. If an excursion of $x_t$ has positive sign, the entrance angle of $\theta_t$ is distributed according to $\mu^+$, and similarly for negative excursions and $\mu^-$. Equivalently, $a$, $\mu^+$, and $\mu^-$ determine a probability measure on $\{-1,1\}\times\Tu$ that gives the entrance behavior of the diffusion from the singularity, but thinking of $\mu^+$ and $\mu^-$ as the conditional distributions of the entrance angle given the sign of the excursion is more consistent with the triangular structure of the system \eqref{Eqn:NaturalScale}. Moreover, since $\theta_t = \theta_0 + \int_0^t \lab x_s\rab^{\alpha} \, dW_s^2$ on the excursion $[0,T_0]$, we see that the $\theta$-process is completely determined by this entrance behavior.

The above gives a complete classification of diffusions on $M_{\alpha}$ (in the present case of $-1<\alpha<0$) that extend Brownian motion on $M$. Indeed, the process is uniquely determined until the first hitting time of the singularity, so such diffusions are determined by their behavior starting from the singularity, and we have the following.

\begin{THM}\label{THM:ConeDiffs}
Let $M_{\alpha}$ be as above for $-1<\alpha<0$, and let $(x_t,\theta_t)$ be a (conservative) diffusion on $M_{\alpha}$ extending Brownian motion on $M$, written in the standard coordinates. Then $(x_t,\theta_t)$ is determined by its behavior starting from $\sing$, which is given by the following parameters: $\gamma\in[0,\infty]$, $a\in[0,1]$, and (Borel) probability measures $\mu^+$ and $\mu^-$ on $\Tu$. More concretely, $x^2_t$ is a squared-Bessel process of dimension $1-\alpha$ on $[0,\infty)$ with reflecting boundary condition at $0$ determined by $\gamma$ (making 0 instantaneously reflecting, slowly reflecting, or absorbing, as above), $x_t$ is recovered from $x^2_t$ by assigning each excursion away from $0$ a positive sign with probability $a$ (and thus a negative sign with probability $1-a$), and on any excursion $t\in(t_1,t_2)$ of $x_t$ away from 0, $\theta_t$ is given as the solution of the SDE $d\theta_t= |x_t|^{\alpha}\, dW^2_t$ with the initial condition $\theta_{t_1}$ distributed as $\mu^+$ if $x_t$ is positive on $(t_1,t_2)$ and $\mu^-$ if $x_t$ is negative on $(t_1,t_2)$. Moreover, if $\gamma=\infty$, $\{x=0\}$ is absorbing and none of the other parameters are relevant, if $\gamma<\infty$ and $a=0$, all excursions of $x_t$ are negative and $\mu^+$ is irrelevant, and if $\gamma<\infty$ and $a=1$, all excursions of $x_t$ are positive and $\mu^-$ is irrelevant, but aside from these exceptions (in what can be considered degenerate cases), there is a one-to-one correspondence between the choice of diffusion and the choice of parameters.
\end{THM}

\subsection{The state space and symmetric extensions}\label{Sect:StateSpace}

The previous classification was restricted to diffusions on $M_{\alpha}$. To clarify, if we start with $M$, since Brownian motion on $M$ explodes toward $\sing=\{x=0\}$ in finite time, to continue the process for all time requires enlarging the state space. In order to make the terminology and relationship to other literature clearer, we call a compactification of $M\cap \{-1\leq x\leq 1\}$ an interior compactification of $M$. The idea is that we want to compactify $M$ at the singularity, but this doesn't, in fact, give a compactification, since $M$ has two more ends, corresponding to $x\rightarrow\pm\infty$. However, we've seen that Brownian motion on $M$ never escapes out of these ends. Thus we want to restrict our attention to a neighborhood of the singularity, and this is what looking at interior compactifications accomplishes.

One way of enlarging the state space is to add a single point for $\{x=0\}$, and this gives the metric space $M_{\alpha}$ that we have been working with, based on the geodesic distance. However, this is not the only possible extension of $M$. Indeed, starting, more functional analytically, from either the Laplacian or the associated Dirichlet form on $M$, one can consider extending the domain of the operator beyond smooth functions compactly supported on $M$. Such extensions are naturally carried by an interior compactification of $M$, but the compactification will depend on the extension and won't necessarily coincide with the one-point compactification that gives $M$. This is the approach followed in \cite{UgoDario}, and we now briefly explain the relationship between their results and the above.

To understand (interior) compactifications of $M$, it is useful to observe that the coordinates $(y,\theta)$ on $M$, where
\[
y=\sgn(x)\frac{1}{\alpha+1}|x|^{\alpha+1}
\]
and $\theta$, of course, is the same $\theta$ from the standard coordinates, give a conformal diffeomorphism from $M$ to the subset of the (Euclidean) cylinder
\[
D= \lp (-\infty,0)\cup(0,\infty)\rp\times \mS^1 \subset \bR\times\mS^1 .
\]
That these coordinates are conformal is already contained in \eqref{Eqn:NaturalScale}, since it shows that Brownian motion on $M$ is a time-change of (Euclidean) Brownian motion on $D$.

The maximal extension of the domain of $\Lap$, as described in \cite{UgoDario}, corresponds to Neumann boundary conditions at the singularity. Unsurprisingly, this corresponds to the ``maximal'' (interior) compactification (from a potential-theoretic viewpoint) of $M$, which is the Martin boundary $\partial_{M}M$ of $M\cap \{-1\leq x\leq 1\}$, or equivalently, of $D\cap \{|y|\leq 1/(\alpha+1)\}$. (Here we put Neumann boundary conditions on $\{x=\pm1\}$ for convenience, in order to restrict the process to $M\cap \{-1\leq x\leq 1\}$.) More concretely, the conformal equivalence with $D$ shows that $\partial_M M$ can be identified with the ``doubled'' Euclidean boundary of $D$, $\{y=0\}$, which is the disjoint union of two copies of $\Tu$.  Thus $\partial_M M$ records both the exit angle $\lim _{t\nearrow T_0}\theta_t$ and exit ``side'' $\lim_{t\nearrow T_0} \sgn(x_t)=\lim_{t\nearrow T_0} \sgn(y_t)$ of the diffusion, where $T_0$ is the first hitting time of $\{y=0\}$ (or the first exit time of $M$) for the diffusion started from a point in $M$. Then the process associated with the Neumann boundary conditions is determined by instantaneous normal reflection of $(y_t,\theta_t)$ back into the component of $D$ the process started in. It's clear how to the construct this diffusion analogously to what was done in the previous section, with $M_{\alpha}$ replaced by $M\cup \partial_M M$, and it's also clear that this process does not descend to a strong Markov process on $M_{\alpha}$.

A second extension of $\Lap$ considered in \cite{UgoDario} is what they call the bridging extension. The corresponding interior compactification is given by identifying pairs of points in $\partial_MM$ with the same $\theta$-coordinate. Since this boundary is also the Euclidean boundary of $D$, we denote it by $\partial_E M$. The corresponding process can be constructed by starting with the diffusion $(|y_t|,\theta_t)$ with Neumann boundary conditions (as just discussed) but assigning signs to the excursions of $(|y_t|,\theta_t)$ randomly with equal probabilities. To see this, note that the condition $\lim_{x\rightarrow 0^+} |x|^{-\alpha}\partial_x u(x,\cdot)=\lim_{x\rightarrow 0^-} |x|^{-\alpha}\partial_x u(x,\cdot) $ in the domain of the Laplacian for the bridging extension becomes $\lim_{y\rightarrow 0^+} \partial_y u(y,\cdot)=\lim_{y\rightarrow 0^-} \partial_y u(y,\cdot) $ after changing coordinates. Then if we consider functions that are even in the $y$ (or $x$) variable, we see that the boundary condition for the $(|y_t|,\theta_t)$-process is just that the normal derivative vanishes, which corresponds to instantaneous normal reflection. Thus the $y_t$ process is independent of the $\theta_t$-process, and we see that the domain of the operator is exactly that of skew-Brownian motion in the trivial case when the skewness vanishes; see Equation (7.6.10) of \cite{ChenFukushimaBook} and the surrounding discussion (that is, the process is just Brownian motion, realized in a slightly non-standard way). This justifies the above claim; alternatively, the process can be thought of as Euclidean Brownian motion on $\overline{D}$ time-changed to spend Lebesgue measure 0 time on $\{y=0\}$ and to solve \eqref{Eqn:NaturalScale} on $D$. Again, the natural state space for this process is $M\cup \partial_E M$, and it does not induce a strong Markov process on $M$.

The third explicit extension of $\Lap$ considered in \cite{UgoDario} is the Friedrich extension, which corresponds to the diffusion killed at the singularity and gives the minimal extension of the domain of $\Lap$. This has $M$ plus a graveyard state for when the particle is killed as its natural state space, which means this case is not covered by Theorem \ref{THM:ConeDiffs}, so the interior compactification of $M$ is ``minimal'' as well, but the diffusion is not conservative. (The case when $\sing$ is absorbing has $\sing$ as a stationary state, so it is not symmetric, and thus not considered in  \cite{UgoDario}.)

Recall that $\omega$ is the Riemannian volume measure on $M$, and let $\ol{\omega}$ be the extension of $\omega$ to $M_{\alpha}$ given by assigning measure 0 to the singularity. Since Brownian motion on a Riemannian manifold is symmetric with respect to the Riemannian volume, it is natural to ask for a diffusion on $M_{\alpha}$ that is symmetric with respect to $\ol{\omega}$. Indeed, both the Neumann and bridging extensions of \cite{UgoDario} are symmetric with respect to the extension of $\omega$ given by assigning measure 0 to $\partial_MM$ or $\partial_EM$. On the other hand, the diffusions of Theorem \ref{THM:ConeDiffs} aren't, in general, symmetric with respect to $\ol{\omega}$. Indeed, the interior compactification taking $M$ into $M_{\alpha}$ gives a one-point compactification of $M\cap \{-1\leq x\leq 1\}$, in the terminology of Chapter 7 of \cite{ChenFukushimaBook} (although what we write as $\omega$ and $\ol{\omega}$ correspond to $\omega_0$ and $\omega$, respectively, in their notation). Thus, according to Theorems 7.5.4 and 7.5.6 of \cite{ChenFukushimaBook}, there is a unique diffusion on $M_{\alpha}$ that extends Brownian motion on $M$ and is symmetric with respect to $\ol{\omega}$. Note that the isometry group of $M_{\alpha}$ is generated by reflection in $x$ and the action of $\mathrm{SO}(2)$ on $\theta$. Then uniqueness implies that such a diffusion, when started from the singularity, must be invariant under $\Iso$. Hence, in Theorem \ref{THM:ConeDiffs}, we must have that $a=1/2$ and both $\mu^+$ and $\mu^-$ are the uniform probability measure on $\mS^1$. Further, since symmetry with respect to $\ol{\omega}$ requires the process to spend 0 time at the singularity, we must have $\gamma=0$. This proves the following.

\begin{THM}
Let $M$ and $\ol{\omega}$ be as above, for $-1<\alpha<0$. Then the unique (conservative) diffusion on $M_{\alpha}$ that extends Brownian motion on $M$, spends time 0 at $\sing$, and is $\ol{\omega}$-symmetric is given by taking $a=1/2$, $\gamma=0$, and $\mu^+$ and $\mu^-$ both to be the uniform probability measure on $\mS^1$ in Theorem \ref{THM:ConeDiffs}. This diffusion is also the unique extension of Brownian motion that spends time 0 at $\sing$ and is invariant under the isometry group of $M_{\alpha}$.
\end{THM}

Let $P_t$ be the semigroup associated to the diffusion described in the preceding theorem. Let $f$ and $g$ be functions in $L^{\infty}(M_{\alpha})$ such that $f=g$ on $M^+$ and for almost every $u<0$,
\[
\int_{\Tu} f(u,\theta)\,d\theta = \int_{\Tu} g(u,\theta)\,d\theta .
\]
Then the $\mathrm{SO(2)}$ invariance of $\mu^-$ means that $P_tf(x)=P_tg(x)$ for every $x>0$ and every $t>0$. (Of course, an analogous result holds for $x<0$ and $u>0$.) In this sense, the diffusion loses information, and only certain ``average'' features of $f$ are communicated across $\sing$. If $\mu^-$ or $\mu^+$ is not uniform, then a similar result holds, except that the $\theta$ averages for each $u$ must be computed with respect to a non-uniform measure (depending on $y$). In any case, the fact that $\sing$ is a single point means that some information must be lost when the process crosses $\sing$.

\subsection{The case $\alpha\leq -1$}\label{Sect:Short}
In the case when $\alpha\leq -1$, the process can enter $M$ from $\sing$, but then never returns (and never hits $\sing$ if it starts from $M$). Thus, if we want a diffusion starting from any point of $M_{\alpha}$, we need only describe how it enters $M$ from $\sing$ (which is a single point). One possibility is for the process to never leave $\sing$, which one can think of as a type of absorbing boundary condition. If the process leaves, it must do so immediately (by the strong Markov property), and just as above, how it enters $M$ is determined by a choice of $a\in[0,1]$ and probability measures $\mu^+$ and $\mu^-$ (unless $a$ is 0  or 1, in which case only one of theses measures is needed).

\section{The case $0\leq \alpha<1$}\label{Sect:Cylinder}

In this case, $M_{\alpha}$ has a cylinder structure at the singularity. In particular, the singularity is now a circle, naturally parametrized by the $\theta$-coordinate, and this would make a complete description of all possible (conservative) diffusions extending Brownian motion on $M$ rather complicated. For instance, such a description is connected to the boundary theory of multidimensional diffusions. More concretely, consider the process $(x^2_t,\theta_t)=(z_t,\theta_t)$ on $[0,\infty)\times \Tu$. Then $z_t$ can undergo sticky, oblique reflection at the boundary, with the parameters determining this reflection depending on $\theta$. Determining a solution without assuming (much) regularity of these parameters (or potentially of the boundary) is a longstanding topic of interest. For example, a construction of a process on a halfspace with general Wentzell boundary conditions was given fairly recently by Watanabe \cite{Watanabe} by extending It\^{o}'s excursion theory, and one can see the references therein for other probabilistic approaches. To extend Brownian motion to $M_{\alpha}$, one would expect a ``two-sided'' version of this type of construction, where the process is potentially sticky at the boundary (and perhaps even diffuses within the boundary) in a way that depends on $\theta$, and when the process re-enters $M$, the distribution of sign of the excursion depends on $\theta$ as does the obliqueness of the ``reflection.'' Constructing such a process, especially for low regularity of the parameters describing this behavior, is well beyond the scope of this note, and it is also in the opposite direction from the more geometrically natural question of determining a ``good'' or ``best'' extension.

Before doing this, motivated by the earlier emphasis on whether or not the process can cross the singularity, we give a simple example to illustrate that the way in which the process crosses the singularity can be unusual. Let $A\subset \Tu$ be a non-empty open subset of $\sing$ such that $A^c$ has non-empty interior. Let $(|x_t|,\theta_t)$ (as a process on $[0,\infty)\times\Tu$) be given by instantaneous normal reflection at the boundary, and let the sign of each excursion of $x_t$ be positive if it begins in $A$ and negative if it begins in $A^c$. Then because the process hits both $A$ and $A^c$ with positive probability from either side of $\sing$, we see that the process will (almost surely) cross $\sing$ infinitely often. However, the crossing is ``non-local,'' in the sense that when the process hits the interior of $A$ from $M^+$, it is distance 0 from $M^-$, but cannot cross into $M^-$ immediately. Instead it must ``go around'' $A$ and cross at $A^c$, and similarly for the process hitting the interior of $A^c$ from $M^-$.

Just as before, the isometry group of $M_{\alpha}$ is $\Iso$. Then we note that if the process spends time 0 at $\sing$ and is symmetric with respect to reflection in $x$, $(|x_t|,\theta_t)$ must be a diffusion on $[0,\infty)\times \Tu$ that reflects instantaneously at the boundary, and $x_t$ can be recovered from $|x_t|$ by giving each excursion a positive or negative sign with probability $1/2$. Additionally, if $(|x_t|,\theta_t)$ is invariant with respect to the $\mathrm{SO}(2)$ action, the reflection must be normal.

Further, the conformal map described in Section \ref{Sect:StateSpace} and given by Equation \eqref{Eqn:NaturalScale} (combined with reflection in $y$) extends to the current case of $0\leq \alpha<1$. Thus we see that the Martin boundary is the same as before, although now it is only ``twice'' the singularity, since $\sing$ is $\partial_EM$ in this case. So while the Neumann extension doesn't give a diffusion on $M$, now the bridging extension does. Also, we see that the construction of the associated diffusion in Section \ref{Sect:StateSpace} (which remains valid here) agrees with the one we just gave under the assumption of invariance under the isometry group. These considerations establish the following.

\begin{THM}
For $0\leq \alpha<1$, the only (conservative) diffusion on $M_{\alpha}$ extending Brownian on $M$ that spends 0 time at $\sing$ and is invariant under the isometry group of $M_{\alpha}$ is given by letting $(|x_t|,\theta_t)$ must be the diffusion on $[0,\infty)\times \Tu$ that undergoes instantaneous normal reflection at the boundary and letting $x_t$ be constructed from $|x_t|$ by giving each excursion a positive or negative sign with probability $1/2$. Moreover, this is the diffusion associated to the bridging extension.
\end{THM}

\def\cprime{$'$}
\providecommand{\bysame}{\leavevmode\hbox to3em{\hrulefill}\thinspace}
\providecommand{\MR}{\relax\ifhmode\unskip\space\fi MR }
\providecommand{\MRhref}[2]{%
  \href{http://www.ams.org/mathscinet-getitem?mr=#1}{#2}
}
\providecommand{\href}[2]{#2}


\begin{thebibliography}{10}

\bibitem{book}
Andrei Agrachev, Davide Barilari, and Ugo Boscain, \emph{A comprehensive
  introduction to sub-riemannian geometry}, Cambridge Studies in Advanced
  Mathematics, vol. 181, Cambridge University Press, Cambridge, 2020.

\bibitem{agrachevboscainsigalotti}
Andrei Agrachev, Ugo Boscain, and Mario Sigalotti, \emph{A
  {G}auss-{B}onnet-like formula on two-dimensional almost-{R}iemannian
  manifolds}, Discrete Contin. Dyn. Syst. \textbf{20} (2008), no.~4, 801--822.
  \MR{2379474 (2009i:53023)}

\bibitem{boscain-laurent}
Ugo Boscain and Camille Laurent, \emph{The {L}aplace-{B}eltrami operator in
  almost-{R}iemannian geometry}, Ann. Inst. Fourier (Grenoble) \textbf{63}
  (2013), no.~5, 1739--1770. \MR{3186507}

\bibitem{UgoDario}
Ugo Boscain and Dario Prandi, \emph{Self-adjoint extensions and stochastic
  completeness of the {L}aplace-{B}eltrami operator on conic and anticonic
  surfaces}, J. Differential Equations \textbf{260} (2016), no.~4, 3234--3269.
  \MR{3434398}

\bibitem{ChenFukushimaBook}
Zhen-Qing Chen and Masatoshi Fukushima, \emph{Symmetric {M}arkov processes,
  time change, and boundary theory}, London Mathematical Society Monographs
  Series, vol.~35, Princeton University Press, Princeton, NJ, 2012.
  \MR{2849840}

\bibitem{ChenFukushimaPaper}
\bysame, \emph{One-point reflection}, Stochastic Process. Appl. \textbf{125}
  (2015), no.~4, 1368--1393. \MR{3310351}

\bibitem{AnjaYor}
Anja G\"oing-Jaeschke and Marc Yor, \emph{A survey and some generalizations of
  {B}essel processes}, Bernoulli \textbf{9} (2003), no.~2, 313--349.
  \MR{1997032}

\bibitem{ItoMcKean}
Kiyoshi It{\^o} and Henry~P. McKean, Jr., \emph{Diffusion processes and their
  sample paths}, Die Grundlehren der Mathematischen Wissenschaften, Band 125,
  Academic Press, Inc., Publishers, New York; Springer-Verlag, Berlin-New York,
  1965. \MR{0199891}

\bibitem{Knight}
Frank~B. Knight, \emph{Essentials of {B}rownian motion and diffusion},
  Mathematical Surveys, vol.~18, American Mathematical Society, Providence,
  R.I., 1981. \MR{613983}

\bibitem{LejaySurvey}
Antoine Lejay, \emph{On the constructions of the skew {B}rownian motion},
  Probab. Surv. \textbf{3} (2006), 413--466. \MR{2280299}

\bibitem{RevuzYor}
Daniel Revuz and Marc Yor, \emph{Continuous martingales and {B}rownian motion},
  Grundlehren der Mathematischen Wissenschaften [Fundamental Principles of
  Mathematical Sciences], vol. 293, Springer-Verlag, Berlin, 1991. \MR{1083357}

\bibitem{Watanabe}
Shinzo Watanabe, \emph{It\^{o}'s theory of excursion point processes and its
  developments}, Stochastic Process. Appl. \textbf{120} (2010), no.~5,
  653--677. \MR{2603058}

\end{thebibliography}
\end{document}